\documentclass[10pt]{article}

\usepackage{theorem,amssymb,amsmath}

\usepackage{graphicx}

\usepackage{hyperref}

\advance \topmargin by -\headheight
\advance \topmargin by -\headsep
\evensidemargin \oddsidemargin
\author{J.-P. Allouche \\
CNRS, IMJ-PRG, Sorbonne \\
4 Place Jussieu \\
F-75252 Paris Cedex 05, France \\
{\tt jean-paul.allouche@imj-prg.fr}
\and
\\
}

\title{H\"older and Kurokawa meet Borwein--Dykshoorn and Adamchik}

\date{ }

\def \proof{\bigbreak\noindent{\it Proof.\ \ }}

\newtheorem{theorem}{Theorem}

\newtheorem{corollary}{Corollary}

\theorembodyfont{\rm}

\newtheorem{remark}{Remark}

\begin{document}

\maketitle

\begin{abstract}
Following our discovery of a nice identity in a recent preprint of Hu and Kim, we show a link 
between the Kurokawa multiple trigonometric functions and two functions introduced respectively 
by Borwein-Dykshoorn and by Adamchik. In particular several identities involving $\zeta(3)$, 
$\pi$ and the Catalan constant $G$ that are proved in these three papers are related.

\medskip


\medskip

\end{abstract}

\section{Introduction}

There is a wealth of special functions arising from geometry and from transcendence theory. 
The purpose of this paper is to provide identities relating some of them. The beginning of 
the story here is a recent preprint of Hu and Kim \cite{HK} that gives the nice identity
$$
\zeta(3) =  \frac{4 \pi^2}{21} \log \left(\frac{e^{\frac{4G}{\pi}} {\mathcal C}_3\left(\frac{1}{4}\right)^{16}}
{\sqrt{2}}\right).
$$
where $G = \sum_{n \geq 0} \frac{(-1)^n}{(2n+1)^2}$ is the Catalan constant,
and ${\mathcal C}_3$ is the Kurokawa-Koyama triple cosine function (see below).

\medskip

The (slightly hidden) occurrences of $7 \zeta(3) / 4 \pi^2$ and of $e^{G/2\pi}$ reminded us 
two (out of four) identities in a paper of Kachi and Tzermias \cite{KT}, namely:
$$
\lim_{n \to \infty} \prod_{k=1}^{2n+1} e^{-1/4} \left(1-\frac{1}{k+1}\right)^{\frac{k(k+1)}{2}(-1)^k}
= \exp\left(\frac{7\zeta(3)}{4\pi^2} + \frac{1}{4}\right) 
$$
and
$$
\lim_{n \to \infty} \prod_{k=1}^{2n+1} \left(1-\frac{2}{2k+1}\right)^{k(-1)^k}
= \exp\left(\frac{2G}{\pi} + \frac{1}{2}\right).
$$
In that paper Kachi and Tzermias proved four identities, and they indicated that they were 
not able to deduce them directly from the values of a function introduced in 1993 
by Borwein and Dykshoorn \cite{BD}. We provided in \cite{Allouche} a proof of their identities, 
using the paper of Borwein and Dykshoorn and their function $D$ for two of the identities, and 
a paper of Adamchik \cite{Adamchik} and his function $E$ for the remaining two. It was thus
tempting to relate these functions $D$ and $E$ to the Kurokawa-Koyama triple cosine.

\medskip

Furthermore, looking at the papers of Kurokawa, we found an expression of $\zeta(3)$
resembling the identity given by Hu et Kim, namely (see \cite[p.~62]{Kurokawa1991},
also see \cite[Theorem~2, p.~209]{zeta3})
$$
\zeta(3) = \frac{8 \pi^2}{7} \log \left(\frac{2^{1/4}}{{\mathcal S}_3(1/2)}\right)
$$
where ${\mathcal S}_3$ is the triple sine function of Kurokawa (in particular 
$$
{\mathcal S}_3(1/2) = e^{1/8} 
\prod_{n \geq 1} \left(\left(1 - \frac{1}{4 n^2}\right)^{n^2} e^{\frac{1}{4}}\right).
$$

It was then even more tempting to study these multiple trigonometric functions and to 
(try to) find a link between all these results.

\bigskip

A possibly surprising fact is that the literature related to the functions of Borwein-Dykshoorn
and of Adamchik appears to be essentially disjoint from the literature related to the Kurokawa
multiple trigonometric functions, with the notable exception of the book \cite{SC} where several
papers of Kurokawa et al. are cited in the references but do not seem to be exploited in the text.
We propose to enlarge the bridge between these two branches of the theory of special functions.

\bigskip

In the sequel, we will recall the definitions of the multiple trigonometric functions first introduced
by Kurokawa in \cite{Kurokawa1991, Kurokawa1992}, then the definitions of the Borwein-Dykshoorn 
function $D$ given in \cite{BD}, and of the Adamchik function $E$ given in \cite{Adamchik}. 
We will obtain close relations between these functions. Furthermore we will show that an identity 
due to Holcombe \cite{Holcombe} can also be obtained using multiple trigonometric functions. 
We will also mention a link between these functions and the generalized Euler constant function 
in \cite{HH} (also see \cite{Xia} and \cite{SH}). We will end with two questions, one of which 
concerning an identity proved in \cite{MS}, that we were not able to address with multiple
trigonometric functions.

\section{Multiple trigonometric functions}

In 1991--1992 (see \cite{Kurokawa1991, Kurokawa1992}) Kurokawa introduced the multiple sine
functions defined by ${\mathcal S}_1(z) := 2 \sin(\pi z)$, and for $r \geq 2$,
$$
{\mathcal S}_r(z) := \exp\left(\frac{z^{r-1}}{r-1} \right) \prod_{n \geq 1} 
\left(P_r\left(\frac{z}{n}\right) P_r\left(-\frac{z}{n}\right)^{(-1)^{r-1}}\right)^{n^{r-1}}
$$
where
$$
P_r(z) = (1-z)\exp\left(z + \frac{z^2}{2} + \ldots + \frac{z^r}{r}\right).
$$
Since $P_r(z) = P_{r-1}(z) e^{z^r/r}$, so that $e^{(z/n)^r} \cdot (e^{(-z/n)^r})^{(-1)^{r-1}} = 1$,
we can clearly write, for $r \geq 2$,
$$
{\mathcal S}_r(z) := \exp\left(\frac{z^{r-1}}{r-1} \right) \prod_{n \geq 1} 
\left( P_{r-1}\left(\frac{z}{n}\right) P_{r-1}\left(-\frac{z}{n}\right)^{(-1)^{r-1}}\right)^{n^{r-1}}.
$$
Note that in particular 
$$
{\mathcal S}_2(z) = e^z \prod_{n \geq 1} \left(\left(\frac{1 - \frac{z}{n}}{1 + \frac{z}{n}}\right)^n e^{2z}\right)
$$ 
is equal to the function $F$ studied by H\"older in \cite[Eq.~(4), p.~515]{Holder}. 
Note that ${\mathcal S}_r$ is equal to $1/\Lambda_r$ for $r > 1$ where 
$\Lambda_r$ is the function defined by Rovinski\u{i} in 1991 (compare 
\cite[(1) p.~62]{Kurokawa1991} and \cite[(5) p.~74]{Rovinskii}). 
Also see \cite{Mellin} and \cite{Rovinsky}.
 
\bigskip

We are not going to give more details in this section about multiple trigonometric functions, except 
for one of the motivations of Kurokawa, who writes in \cite{Kurokawa1991}: {\it as an application we 
report the calculation of the "gamma factors" of Selberg-Gangolli-Wakayama zeta functions of rank 
one locally symmetric spaces}, and for pointers to a limited number of references, e.g., 
\cite{Vigneras} for an early approach of a related question, and the survey of Manin \cite{Manin}.

\section{The Borwein-Dykshoorn function}

The Borwein-Dykshoorn function
$$
D(x) = \lim_{n \to \infty} \prod_{k=1}^{2n+1} \left(1+\frac{x}{k}\right)^{k(-1)^{k+1}}
$$
was introduced in \cite{BD} as a generalization of a result of Melzak \cite{Melzak}
proving that
$$
\lim_{n \to \infty} \prod_{k=1}^{2n+1} \left(1+\frac{2}{k}\right)^{k(-1)^{k+1}} = \frac{\pi e}{2}\cdot
$$
This function can be extended to a meromorphic function on ${\mathbb C}$ with poles at the
negative even integers. We prove that this function is related to ${\mathcal S}_2$.

\begin{theorem}\label{functionD}
The following equality holds:
$$
\frac{D(x)}{D(-x)} = e^x \ \frac{{\mathcal S}_2(x/2)^4}{{\mathcal S}_2(x)}\cdot
$$

\end{theorem}

\proof Let $D_n(x)$ be defined by
$$
D_n(x) :=  \prod_{k=1}^{2n+1} \left(1+\frac{x}{k}\right)^{k(-1)^{k+1}}\cdot
$$
We can write
$$
D_n(x) = \prod_{j=1}^n \left(1+\frac{x}{2j}\right)^{-2j} 
\prod_{j=0}^n \left(1+\frac{x}{2j+1}\right)^{2j+1}
= \prod_{j=1}^n \left(1+\frac{x}{2j}\right)^{-2j} 
\frac{\displaystyle\prod_{k=1}^{2n+1} \left(1+\dfrac{x}{k}\right)^k}
{\displaystyle\prod_{k=1}^n \left(1+\dfrac{x}{2k}\right)^{2k}}\cdot
$$
Thus
$$
D_n(x) = \frac{\displaystyle\prod_{k=1}^{2n+1}\left(1+\dfrac{x}{k}\right)^k} 
{\displaystyle\prod_{j=1}^n \left(1+\dfrac{x}{2j}\right)^{4j}}
$$
which implies
$$
\frac{D_n(x)}{D_n(-x)} = \prod_{k=1}^{2n+1}\left(\frac{1+\dfrac{x}{k}}{1-\dfrac{x}{k}}\right)^k
\left(\prod_{k=1}^n\left(\frac{1-\dfrac{x/2}{k}}{1+\dfrac{x/2}{k}}\right)^k\right)^4\cdot
$$
Multiplying the first product by $e^{-2x(2n+1)}$ and the second product by
$(e^{\frac{x}{2}} e^{nx})^4$ does not change the quantity $D_n(x)D_n(-x)$, thus
$$
\frac{D_n(x)}{D_n(-x)} = 
\prod_{k=1}^{2n+1}\left(e^{-2x}\left(\frac{1+\dfrac{x}{k}}{1-\dfrac{x}{k}}\right)^k\right)
\left(e^{\frac{x}{2}}\prod_{k=1}^n \left(e^x \left(\frac{1-\dfrac{x/2}{k}}{1+\dfrac{x/2}{k}}\right)^k\right)\right)^4
$$
which gives, when $n$ tends to infinity,
$$
\frac{D(x)}{D(-x)} = e^x \ \displaystyle\dfrac{{\mathcal S}_2(x/2)^4}{{\mathcal S}_2(x)}\cdot \ \ \ \ \ \Box
$$

A corollary of this result gives the value of an infinite product studied in \cite{KT} which was proved
again in \cite{Allouche} by using the Borwein--Dykshoorn function.

\begin{corollary}
The following identities hold
$$
\lim_{n \to \infty} \prod_{k=1}^{2n} \left(1-\frac{2}{2k+1}\right)^{k(-1)^k}
= \ \exp\left(\frac{2G}{\pi} - \frac{1}{2}\right)
$$
and
$$
\lim_{n \to \infty} \prod_{k=1}^{2n+1} \left(1-\frac{2}{2k+1}\right)^{k(-1)^k}
= \ \exp\left(\frac{2G}{\pi} + \frac{1}{2}\right)
$$
where $G = \displaystyle\sum_{n \geq 0} \frac{(-1)^n}{(2n+1)^2}$ is the Catalan constant.
\end{corollary}

\proof  Since $\left(1 -  \dfrac{2}{4n+3}\right)^{-(2n+1)}$ tends to $e$ (take the logarithm), it
suffices to prove the second identity. We write
$$
\prod_{k=1}^{2n+1} \left(1-\frac{2}{2k+1}\right)^{k(-1)^k} = \ \ 
\prod_{k=1}^{2n+1} \left(\frac{2k+1}{2k-1}\right)^{k(-1)^{k+1}}
= \ \ \frac{\displaystyle\prod_{k=1}^{2n+1}\left(1 + \dfrac{1}{2k}\right)^{k(-1)^{k+1}}}
    {\displaystyle\prod_{k=1}^{2n+1}\left(1 - \dfrac{1}{2k}\right)^{k(-1)^{k+1}}}\cdot
$$
Hence, using Theorem~\ref{functionD}
$$
\lim_{n \to \infty} \prod_{k=1}^{2n+1} \left(1-\frac{2}{2k+1}\right)^{k(-1)^k} = \ \ \frac{D(1/2)}{D(-1/2)}
\ = \ e^{1/2} \frac{{\mathcal S}_2(1/4)^4}{{\mathcal S}_2(1/2)}\cdot
$$
But we have ${\mathcal S}_2(1/2) = 2^{1/2}$ and ${\mathcal S}_2(1/4)
= 2^{1/8}\exp\left(\frac{G}{2\pi}\right)$ (see \cite[p.~852]{KK2003}, hence
$$
\lim_{n \to \infty} \prod_{k=1}^{2n+1} \left(1-\frac{2}{2k+1}\right)^{k(-1)^k} = \ \
e^{1/2} \frac{{\mathcal S}_2(1/4)^4}{{\mathcal S}_2(1/2)} \ = \ 
\exp\left(\frac{2G}{\pi} + \frac{1}{2}\right)\cdot
$$

\begin{remark}
\ { }
\begin{itemize}

\item The function $D$ can also be written (see, e.g., \cite{Allouche})
$$
D(x) = e^x\prod_{k \geq 1}\left(e^{-x} \left(1 + \dfrac{x}{k}\right)^k\right)^{(-1)^{k+1}}
$$
 \item
 The left hand term of the equality in Theorem~\ref{functionD}, say $f(x) := D(x)/D(-x)$, has the
property that $f(-x) = 1/f(x)$. Hence this is the same for the right hand term, but it is of course 
easy to see that the function ${\mathcal S}_2(x)$ itself satisfies the identity
${\mathcal S}_2(-x) = ({\mathcal S}_2(x))^{-1}$ (see, e.g., \cite[Eq.(5), p.~516]{Holder}).

\item 
Note that \cite[Theorem, p.~204]{BD} gives an explicit (finite) formula in terms of the function 
$\Gamma$ and of the generalized Gamma function $\Gamma_1$ defined in \cite{Bendersky}
and itself closely related to the Barnes function. There is also a close formula for $D(x)$ in 
\cite[Proposition~5, p.~284]{Adamchik}. Also see \cite[Example~20, p.~139]{MS}. Finally, 
an expression of $D(x)$ in terms of the function $\gamma_{\alpha}$ and its derivative is 
given in \cite[p.~86]{Allouche}, where $\gamma_{\alpha}$ (see \cite{HH} and \cite{SH}; 
also see \cite{Xia}) is defined by:
$$
\gamma_{\alpha}(z) := \sum_{n \geq 1} z^{n-1}
\left(\alpha - n \log \left(1 + \frac{\alpha}{n}\right)\right).
$$

\end{itemize}
\end{remark}

\section{The Adamchik function $E$}

Adamchik defined in \cite{Adamchik} a function $E$ by
$$
E(x) := \lim_{n \to \infty} \prod_{k=1}^{2N} \left( 1 - \dfrac{4x^2}{k^2} \right)^{k^2(-1)^{k+1}}.
$$
It is easy to see that $E(x)$ has also the following expression:
$$
E(x) = \prod_{k \geq 1} \left(e^{4x^2}\left(1 - \dfrac{4x^2}{k^2}\right)^{k^2}\right)^{(-1)^{k+1}}.
$$
Recall that ${\mathcal S}_3(x)$ is defined by
$$
{\mathcal S}_3(x) = 
e^{\frac{x^2}{2}} \prod_{n \geq 1} \left(e^{x^2}\left(1 - \frac{x^2}{n^2}\right)^{n^2}\right).
$$
Also recall that the function ${\mathcal C}_3(x)$ is defined by (see, e.g., \cite[p.~1]{KW2003})
$$
{\mathcal C}_3(x) = \prod_{n \geq 1} 
\left(e^{x^2} \left(1 - \dfrac{4x^2}{(2n- 1)^2}\right)^{\frac{(2n-1)^2}{4}}\right)
$$
Now we state the following result.
\begin{theorem}\label{adamchik-E}
We have the following identity:
$$
E(x) = e^{2x^2} \dfrac{{\mathcal S}_3(2x)}{{\mathcal S}_3(x)^8} 
= e^{2x^2} \dfrac{{\mathcal C}_3(x)^8}{{\mathcal S}_3(2x)}\cdot 
$$
\end{theorem}

\proof
We write
$$
E(x) = \prod_{k \geq 1} \left(e^{4x^2}\left(1 - \dfrac{4x^2}{k^2}\right)^{k^2}\right)^{(-1)^{k+1}}
= \ \ \frac{\displaystyle\prod_{\stackrel{\scriptstyle k \geq 1}{k \ \text{odd}}} 
\left(e^{4x^2}\left(1 - \dfrac{4x^2}{k^2}\right)^{k^2}\right)}
{\displaystyle\prod_{k \geq 1} 
\left(e^{4x^2}\left(1 - \dfrac{x^2}{k^2}\right)^{4k^2}\right)} \ \ 
= \ \ \frac{\displaystyle\prod_{k \geq 1} \left(e^{4x^2}\left(1 - \dfrac{4x^2}{k^2}\right)^{k^2}\right)}
{\displaystyle\prod_{k \geq 1} 
\left(e^{4x^2}\left(1 - \dfrac{x^2}{k^2}\right)^{4k^2}\right)^2}\cdot
$$
Hence
$$
E(x) = \left({\mathcal S}_3(2x) e^{-2x^2}\right) \left({\mathcal S}_3(x) e^{-x^2/2}\right)^{-8}
= e^{2x^2} \frac{{\mathcal S}_3(2x)}{{\mathcal S}_3(x)^8}\cdot
$$
Now (see, e.g., \cite[Theorem~1.1, p.~123]{KW2003}),
$$
{\mathcal S}_3(2x) = {\mathcal S}_3(x)^4 {\mathcal C}_3(x)^4,
$$
which gives
$$ 
E(x) = e^{2x^2} \frac{{\mathcal C}_3(x)^8}{{\mathcal S}_3(2x)} \ \ \ \ \ \Box
$$

\begin{remark}\label{KK-adamchik}

\ { }
\begin{itemize}
\item[$*$]
The first equality in Theorem~\ref{adamchik-E} above and Equation~(2.15) in \cite[p.~850]{KK2003}
imply Corollary~1 in the paper of Adamchik \cite[p.~286]{Adamchik} where $E(x)$ is expressed
in terms of ${\text Li}_2(\pm e^{2i\pi x})$ and ${\text Li}_3(\pm e^{2i\pi x})$, up to replacing $i$ 
with $-i$, which does not change the result since $E(x)$ is real
(recall that Li$_k(x) := \sum_{n \geq 1} \frac{\sin nx}{n^k}$). It is also interesting to compare
these expressions and \cite[Proposition~2 p.~618]{Zagier}.
\item[$*$]
Note that \cite[Proposition~6, p.~286]{Adamchik} gives a close formula for $E$ in terms of the 
Barnes function and of the generalized gamma function $\Gamma_3$.
\end{itemize}
\end{remark}

We have as a corollary the following particular cases given in \cite[p.~287]{Adamchik}.

\begin{corollary}
We have
$$
\ \ \ \lim_{N \to \infty} \prod_{k=1}^{2N} \left(1 - \dfrac{1}{4k^2}\right)^{k^2(-1)^{k+1}}
= \ \ \exp \left(\frac{1}{8} - \frac{2G}{\pi} + \frac{7 \zeta(3)}{2\pi^2}\right)
$$
and
$$
\lim_{N \to \infty} \prod_{k=2}^{2N} \left(1 - \dfrac{1}{k^2}\right)^{k^2(-1)^{k+1}}
= \ \ \frac{\pi}{4} \exp \left(\frac{1}{2} + \frac{7 \zeta(3)}{\pi^2}\right)
$$ 
\end{corollary}

\proof   

The left hand side of first equality is equal to $E(1/4)$ thus, using 
Theorem~\ref{adamchik-E}, to $e^{1/8} {\mathcal C}_3(1/4)^8 {\mathcal S}_3(1/2)^{-1}$.
But ${\mathcal S}_3(1/2)$ is given in \cite[p.~206]{zeta3}, where we find
$$
\zeta(3) = \frac{8\pi^2}{7} \log(2^{1/4} A^{-1}), \ \text{where} \ 
A := e^{1/8} \prod_{n \geq 1}\left( e^{1/4} \left(1 - \dfrac{1}{4n^2}\right)^{n^2}\right)
$$
which can be written 
$$
{\mathcal S}_3(1/2) = 2^{1/4} \exp\left(-7\dfrac{\zeta(3)}{8\pi^2}\right).
$$
The value of ${\mathcal C}_3(1/4)$ is given in \cite[p.~1]{KW2003}, where we find
$$
{\mathcal C}_3(1/4) = 2^{1/32} \exp\left(\dfrac{21\zeta(3)}{64 \pi^2} - \dfrac{G}{4\pi}\right).
$$
Here $G$ is the Catalan constant. Note that there seems to be a misprint in \cite{KW2003},
namely $G/4\pi$ there should be replaced with $-G/4\pi$ as above (also see \cite{HK}). 
So we finally have
$$
E(1/4) = \exp\left(\dfrac{1}{8} - \dfrac{2G}{\pi} + \dfrac{7\zeta(3)}{2\pi^2}\right).
$$

\section{Holcombe's infinite product}

In a 2013 paper \cite{Holcombe} Holcombe proved that
$$
\pi = e^{3/2} \prod_{n \geq 2} \left(e\left(1 - \frac{1}{n^2}\right)^{n^2}\right).
$$
This result can be deduced from the value of the derivative at $1$ of the triple trigonometric
function ${\mathcal S}_3$. Namely
$$
{\mathcal S}_3(x) = 
e^{\frac{x^2}{2}} \prod_{n \geq 1} \left(e^{x^2}\left(1 - \frac{x^2}{n^2}\right)^{n^2}\right)
$$
so that ${\mathcal S}_3(1) = 0$ and
$$
\frac{{\mathcal S}_3(x) - {\mathcal S}_3(1)}{x-1} =
- e^{\frac{3}{2}} (1+x)\prod_{n \geq 2} \left(e\left(1 - \frac{1}{n^2}\right)^{n^2}\right).
$$
Letting $x$ tend to $1$, we thus have
$$
{\mathcal S}_3'(1) =  
- 2  e^{\frac{3}{2}} \prod_{n \geq 2} \left(e\left(1 - \frac{1}{n^2}\right)^{n^2}\right).
$$
Using the value ${\mathcal S}_3'(1) = -2\pi$ 
(see \cite[Prop.~4.1, p.~125]{KW2003}) gives Holcombe's result.

\section{Miscellanea}

Further links between the above mentioned functions and other special functions can
be found in the literature. We give only quick remarks in this section.

\bigskip

Recall that the dilogarithm function, already alluded to in Remark~\ref{KK-adamchik},
and the Clausen function are respectively defined by
$$
\text{Li}_2(z) := \sum_{n \geq 1} \frac{z^n}{n^2} \ \ (\text{where} \ |z| \leq 1),
\ \text{and} \ \ 
\text{Cl}_2(\theta) = \sum_{n \geq 1} \frac{\sin n \theta}{n^2}\cdot
$$
The Clausen function Cl$_2$ is related to the Barnes function $G$ (see, e.g.,
\cite[Equation~(2) p.~175]{SC}):
$$
\text{Cl}_2(x) = x \log \pi - x \log\left(\sin \left(\frac{x}{2}\right)\right)
+ 2\pi \log \frac {G(1 - \frac{x}{2\pi})}{G(1 + \frac{x}{2\pi})}\cdot
$$
But we can find in \cite{KK2003} the relation
$$
{\mathcal S}_2(z) = 
(2 \sin \pi z)^z \exp\left(\frac{1}{2\pi} \sum_{n \geq 1} \frac{\sin(2\pi n z)}{n^2}\right)
$$
hence
$$
{\mathcal S}_2(z) = (2 \sin \pi z)^z \exp\left(\frac{1}{2\pi} \text{Cl}_2(2\pi z)\right).
$$
Thus\
$$
\hskip 3truecm \log\frac{G(1+t)}{G(1-t)} = t \log(2\pi) - \log{\mathcal S}_2(t). \hskip 3truecm (*)
$$
This identity can also be obtained directly from a relation for ${\mathcal S}_2(x)$ given in
\cite[Theorem~2.5, p.~847]{KK2003} and a relation given in \cite[(2.1), p.~94]{CS}, namely
$$
{\mathcal S}_2(x) = \exp\left(\int_0^x \pi t \cot(\pi t) \ dt\right)
\ \
\text{and}
\ \ 
\int_0^x \pi t \cot (\pi t) dt = x \log(2\pi) + \log \frac{G(1-x)}{G(1+x)}\cdot
$$

\medskip 

Relation $(*)$ can be double checked with the values of ${\mathcal S}_2(1/2) = \sqrt{2}$ and
${\mathcal S}_2((1/4)= 2^{1/8} e^{G/2\pi}$ given in \cite[Examples~2.9, p.~852]{KK2003} and 
the values of $G(3/2)/G(1/2) = \Gamma(1/2) = \sqrt{\pi}$ and 
$G(5/4)/G(3/4) = 2^{1/8} \pi^{1/4} e^{-G/2\pi}$ 
given in \cite[p.~271]{Markov}. This relation can also be used to give a closed form of values 
of ${\mathcal S}_2(t)$ using known values of $G(1+t)/G(1-t)$, e.g., for $t = 1/8$ or $t = 1/12$
(see \cite[Section~5]{CSA} and the corrections and results in \cite[p.~272--273]{Markov}).
 
\begin{remark}

In \cite[p.~94]{CS} the authors write that the relation 
$$
\int_0^x \pi t \cot (\pi t) dt = x \log(2\pi) + \log \frac{G(1-x)}{G(1+x)}\cdot
$$
is originally due to Kinkelin. Actually Kinkelin proves a similar result ---see the first 
equality after (26.) p.~135 in \cite{Kinkelin}--- for a function which is actually equal 
to $\Gamma(x)^{x-1} / G(x)$: as written in \cite[p.~136]{Finch} this function {\it has 
been comparatively neglected by researchers in favor of $G$.}

\medskip

This discussion opens the way to other links with several other special functions (e.g., the
dilogarithm already mentioned, the inverse tangent integral 
Ti$_2(z) := \sum_{j \geq 1} (-1)^{j+1} y^j/j^2$ and the Legendre 
chi-function $\chi_2(z) := \sum_{k \geq 0} y^{2k+1}/(2k+1)^2$, both defined for $|z| \leq 1$)
but we will not go further in this direction: the reader can consult in particular \cite{CSA, Markov} 
and the references therein, with a special mention for the book of Lewin \cite{Lewin}.

\end{remark}

\section{Conclusion}

We end this paper with two questions.

\medskip

First, we note that combining the expression of $D(x)/D(-x)$ given above in terms of 
${\mathcal S}_2(x/2)$ and ${\mathcal S}_2(x)$ and the expression of $D(x)$ in terms
of the Hessami Pilehrood-Sondow-Hadjicostas function given in \cite{Allouche}, one 
can obtain a relation between this function and ${\mathcal S}_2(x)$. Is it possible to
obtain a relation with the higher multiple trigonometric functions ${\mathcal S}_r(x)$?

\medskip

Second, an interesting identity is given in \cite[p.~125]{MS} where it is proved by unusual 
methods:
$$
\lim_{n \to \infty} \left[e^{\frac{n}{4}(4n+1)} n^{-\frac{1}{8} - n(n+1)} (2\pi)^{-\frac{n}{2}}
\prod_{k=1}^{2n} \Gamma\left(1 + \dfrac{k}{2}\right)^{k(-1)^k}\right] = 
2^{\frac{1}{12}} \exp \left(\dfrac{5}{24} - \dfrac{3}{2} \zeta'(-1) - \dfrac{7 \zeta(3)}{16 \pi^2}\right)
$$
where the right hand constant can also be written
$$
(2e)^{1/12} A^{3/2} \exp\left(- \frac{7 \zeta(3)}{16 \pi^2}\right)
$$
with $A := \exp(\frac{1}{12} - \zeta'(-1))$ is the Glaisher-Kinkelin constant. We did not succeed 
in finding a relation between this identity with multiple trigonometric functions. Of course the
quantity $7\zeta(3)/4\pi^2$ occurs, but, after all, this is nothing but $-7 \zeta'(-2)$ which looks
somehow more mundane. Is it possible to find a proof of this identity that uses multiple
trigonometric functions?

\bigskip

\noindent
{\bf Acknowledgments} \ \ We would like to thank S. Hu and M.-S. Kim for discussions about
their preprint \cite{HK}. We would also like to thank Olivier Ramar\'e and Sanoli Gun, in 
particular for their help in obtaining some of the references.










\end{document}